\documentclass[12pt,a4paper]{article}
\usepackage{amssymb,amsfonts,amsmath,times}

\newtheorem {theorem}{Theorem}[section]
\newtheorem {corollary}[theorem]{Corollary}
\newtheorem {lemma}[theorem]{Lemma}

\newtheorem {definition}[theorem]{Definition}
\newtheorem {remark}[theorem]{Remark}
\newtheorem {condition}[theorem]{Condition}
\newtheorem {example}[theorem]{Example}
\numberwithin {equation}{section}

  \allowdisplaybreaks [3]

\begin{document}

\title{Stochastic differential equations in a scale of Hilbert spaces}
\author{Alexei Daletskii \\
Department of Mathematics, University of York, UK}

\maketitle

\begin{abstract}
A stochastic differential equation with coefficients defined in a scale of
Hilbert spaces is considered. The existence and uniqueness of finite time
solutions is proved by an extension of the Ovsyannikov method. This result
is applied to a system of equations describing non-equilibrium stochastic
dynamics of (real-valued) spins of an infinite particle system on a typical
realization of a Poisson or Gibbs point process in ${\mathbb{R}}^{n}$.
\end{abstract}

\def\MSC{\par\leavevmode\hbox {\it MSC2010:\ }}%
\def\keywords{\par\leavevmode\hbox {\it Keywords:\ }}%





\section{Introduction}

Evolution differential and stochastic differential equations in Banach
spaces play hugely important role in many parts of mathematics and its
applications. This class of equations unifies infinite systems of ordinary
differential equations and partial differential equations (realized in $%
l_{p} $-type spaces of sequences and Sobolev-type spaces, respectively), and
their stochastic counterparts, see e.g. \cite{Deim}, \cite{DZ} and
references therein and modern developments in e.g. \cite{DDFR}.

So let us consider a stochastic differential equation (SDE) of the form
\begin{equation}
d\xi (t)=f(\xi (t))dt+B(\xi (t))dW(t)  \label{sde}
\end{equation}%
in a Banach space $X$, where $f$ and $B$ are given vector and operator
fields on $X$ respectively and $W$ a suitable Wiener process in $X$. The
standard approach to such equations usually requires that $f=A+\phi $, where
(C1) $A$ is a generator of a $C_{0}$-semigroup in $X$, and (C2) $\phi $ and $%
B$ satisfy certain Lipschitz or dissipativity conditions in $X$. Then the
existence, uniqueness and regularity of solutions of the corresponding
Cauchy problem can be proved.

This classical theory does not cover some important examples motivated by
e.g. problems of statistical mechanics and hydrodynamics. In particular,
there are situations where $A$ fails to satisfy condition (C1) but is
instead bounded in a scale of Banach spaces $X_{\alpha }$, $\alpha \in {%
\mathcal{A}}$, where ${\mathcal{A\subset }}{\mathbb{R}}^{1}$ is an interval
and $X_{\alpha }\subset X_{\beta }$ if $\alpha \leq \beta $. That is, $A$ is
a bounded operator acting from $X_{\alpha }$ to $X_{\beta }$ for any $\alpha
<\beta $, and
\begin{equation}
\left \Vert Ax\right \Vert _{X_{\beta }}\leq c\left (\beta -\alpha \right
)^{-1}\left \Vert x\right \Vert _{X_{\alpha }}  \label{A-cond}
\end{equation}
for all $x\in X_{\alpha }$ and some constant $c>0$ (independent of $\alpha $
and $\beta $ but possibly dependent on the interval ${\mathcal{A}}$).

In this framework, equation (\ref{sde}) with no diffusion term ($B\equiv 0$)
has been studied by Ovsyannikov's method, see e.g. \cite{Deim} and modern
developments and references in \cite{Fink}, \cite{BHP}. Moreover, instead of
(C2), the non-linear drift term $\phi $ is allowed to satisfy a generalized
Lipschitz condition in the scale $(X_{\alpha })_{\alpha \in {\mathcal{A}}}$
with singularity of the type as in (\ref{A-cond}) (see \cite{Ni,Saf,BHP}).
The price to pay here is that the existence of a solution with initial value
in $X_{\alpha }$ can only be proved in the bigger space $X_{\beta }$, $\beta
>\alpha $. The lifetime of this solution depends on $\alpha $ and $\beta $
(and the interval ${\mathcal{A}}$ itself).

The aim of the present work is to extend Ovsyannikov's method to the case of
stochastic differential equations. We require the drift $f$ to be a map from
$X_{\alpha }$ to $X_{\beta }$ for any $\alpha <\beta $ and satisfy a
generalized Lipschitz condition with singularity $\left (\beta -\alpha
\right )^{-1/2}$ (and make similar assumption about the diffusion
coefficient $B$), see Condition \ref{Lip-drift} given in the next section,
and prove the existence and uniqueness of finite time solutions of the
corresponding Cauchy problem. Observe that the singularity allowed here is
weaker than in the deterministic case (cf. (\ref{A-cond})), which is related
to the specifics of the Ito integral estimates. As in the deterministic
case, the solution will live in the scale $X_{\alpha }$, $\alpha \in {%
\mathcal{A}}$. For simplicity, we assume that all $X_{\alpha }$ are Hilbert
spaces, although all our results hold in a more general situation of
suitable Banach spaces. The proof is based on the contractivity of the
corresponding integral transformation of a weighted space of trajectories in
$\cup _{\alpha \in {\mathcal{A}}}X_{\alpha }$ (constructed similar to the
ones used in \cite{Ni,Saf,BHP}).

Our main example is motivated by the study of countable systems of particles
randomly distributed in a Euclidean space ${\mathbb{R}}^{n}$ (of the type
considered in \cite{DKKP}, \cite{DKKP1}). Each particle is characterized by
its position $x$ and an internal parameter (spin) $\sigma _{x}\in {\mathbb{R}%
}$. For a given fixed (\textquotedblleft quenched\textquotedblright )
configuration $\gamma $ of particle positions, which is a locally finite
subset of ${\mathbb{R}}^{n}$, we consider a system of stochastic
differential equations describing (non-equilibrium) dynamics of spins $%
\sigma _{x},$ $x\in \gamma $. Two spins $\sigma _{x}$ and $\sigma _{y}$ are
allowed to interact via a pair potential if the distance between $x$ and $y$
is no more than a fixed interaction radius $r$, that is, they are neighbors
in the geometric graph defined by $\gamma $ and $r.$ Vertex degrees of this
graph are typically unbounded, which implies that the coefficients of the
corresponding equations cannot be controlled in a single Hilbert or Banach
space (in contrast to spin systems on a regular lattice, which have been
well-studied, see e.g. \cite{DZ1} and modern developments in \cite{INZ}, and
references therein). However, under mild conditions on the density of $%
\gamma $ (holding for e.g. Poisson and Gibbs point processes in ${\mathbb{R}}%
^{n}$), it is possible to apply the approach discussed above and construct a
solution in the scale of Hilbert spaces $S_{\alpha}^{\gamma }$ of weighted
sequences $(q_x)_{x\in\gamma }$ such that $\sum_{x\in \gamma }\left\vert
q_{x}\right\vert ^{2}e^{-\alpha \left\vert x\right\vert }<\infty ,\ \alpha
>0 $.

Observe that the family $X_\alpha=S^\gamma_\alpha$, $\alpha>0$, forms the
dual to nuclear space $\Phi ^{\prime }=\cup _{\alpha }X_\alpha$. SDEs
 on such spaces were considered in \cite{KMW}, \cite{KJX}. The
existence of solutions to the corresponding martingale problem was proved
under assumption of continuity of coefficients on $\Phi ^{\prime }$ and
their linear growth (which, for the diffusion coefficient, is supposed to
hold in each $\alpha$-norm). Moreover, the existence of strong solutions
requires a dissipativity-type estimate in each $\alpha$-norm, too, which
does not hold in our framework.

In the last subsection, we prove the uniqueness of the infinite-particle
dynamics using more classical methods, which generalise those applied to
deterministic systems in \cite{LLL}, \cite{DaF}.

\bigskip

\section{Setting}

Let us consider a family of separable Hilbert spaces $X_{\alpha }$ indexed
by $\alpha \in \left [\alpha_{\ast},\alpha^{\ast}\right ]$ with fixed $0\le
\alpha_{\ast},\alpha^{\ast}<\infty $, and denote by $\left \Vert \cdot
\right \Vert _{\alpha }$ the corresponding norms. We assume that
\begin{equation}
X_{\alpha }\subset X_{\beta }\ {\text {and }}\left \Vert u\right \Vert
_{\beta }\leq \left \Vert u\right \Vert _{\alpha }{\text { if }}\alpha
<\beta ,\ u\in X_{\beta },  \label{scale}
\end{equation}
where the embedding means that $X_{\alpha }$ is a vector subspace of $%
X_{\beta }$. When speaking of these spaces and related objects, we will
always assume that the range of indices is $\left
[\alpha_{\ast},\alpha^{%
\ast}\right ]$, unless stated otherwise.

Let $W(t)$ be a cylinder Wiener process in a separable Hilbert space ${%
\mathcal{H}}$ defined on a suitable filtered probability space. Introduce
notation
\begin{equation*}
H_{\beta }\equiv HS({\mathcal{H}},X_{\beta }):=\left\{ {\text{%
Hilbert-Schmidt operators }}{\mathcal{H}}\rightarrow X_{\beta }\right\} .
\end{equation*}%
We will denote by $\left\Vert \cdot \right\Vert _{H_{\beta }}$ its standard
norm. Our aim is to construct a strong solution of equation (\ref{sde}),
that is, a solution of the stochastic integral equation
\begin{equation}
u(t)=u_{0}+\int_{0}^{t}f(u(s))ds+\int_{0}^{t}B(u(s))dW(s),  \label{SDE-int}
\end{equation}%
with coefficients acting in the scale of spaces (\ref{scale}). More
precisely, we assume that $f:X_{\alpha }\rightarrow X_{\beta }$ and $%
B:X_{\alpha }\rightarrow H_{\beta }$ for any $\alpha <\beta $, and the
following Lipschitz-type condition is satisfied.

\begin{condition}
\label{Lip-drift}There exists a constant $L$ such that
\begin{equation}
\left\Vert f(u)-f(v)\right\Vert _{\beta }\leq \frac{L}{\left\vert \beta
-\alpha \right\vert ^{1/2}}\left\Vert u-v\right\Vert _{\alpha }  \label{lg}
\end{equation}%
and
\begin{equation}
\left\Vert B(u)-B(v)\right\Vert _{H_{\beta }}\leq \frac{L}{\left\vert \beta
-\alpha \right\vert ^{1/2}}\left\Vert u-v\right\Vert _{\alpha }
\label{lg-diff}
\end{equation}%
for any $\alpha <\beta $ and all $u,v\in X_{\alpha }$.
\end{condition}

We denote by ${\mathcal{GL}}^{(1)}$ and ${\mathcal{GL}}^{(2)}$ the sets of
mappings $f$ and $B$ under conditions (\ref{lg}) and (\ref{lg-diff}),
respectively.

\begin{remark}
The Lipschitz constant $L$ may depend on $\alpha ^{\ast }$ and $\alpha
_{\ast }$, as usually happens in applications.
\end{remark}

\begin{remark}
In contrast to the classical Ovsyannikov method for deterministic equations,
where the right-hand side of (\ref{lg}) is proportional to $\left (\beta
-\alpha \right )^{-1}$, we have to require stronger bounds with the
singularity $\left (\beta -\alpha \right )^{-1/2}$. This is due to the
presence of the Ito stochastic integral in (\ref{SDE-int}).
\end{remark}

\begin{remark}
\label{rem-bound}Setting $v=0$ in (\ref{lg}) and (\ref{lg-diff}), we obtain
linear growth conditions
\begin{equation*}
\left\Vert f(u)\right\Vert _{\beta }\leq \frac{K}{\left\vert \beta -\alpha
\right\vert ^{1/2}}\left( 1+\left\Vert u\right\Vert _{\alpha }\right)
\end{equation*}%
and
\begin{equation*}
\left\Vert B(u)\right\Vert _{H_{\beta }}\leq \frac{K}{\left\vert \beta
-\alpha \right\vert ^{1/2}}\left( 1+\left\Vert u\right\Vert _{\alpha }\right)
\end{equation*}%
for some constant $K$, any $\alpha <\beta $ and all $u\in X_{\alpha }$.
\end{remark}
 
\begin{remark}
Assume that $\phi $ is Lipschitz continuous in each $X_{\alpha }$ with a
uniform Lipschitz constant $M$. Then $\phi \in {\mathcal{GL}}^{(1)}$ with $L=%
\sqrt{\alpha ^{\ast }-\alpha _{\ast }}M$.
\end{remark}

\begin{remark}
\label{rem-scale}Some authors have used the scale $X_{\alpha }$ such that $%
X_{\alpha }\subset X_{\beta }$ if $\alpha >\beta $. This framework can be
transformed to our setting by an appropriate change of the parametrization,
e.g. $\alpha \mapsto \alpha^{\ast}-\alpha $.
\end{remark}

\section{Main results\label{sec-main}}

Let us fix $b>0$ and define the function
\begin{equation*}
p_{b}(\alpha ,t):=1-\left (\left (\alpha -\alpha_{\ast}\right )b\right
)^{-1}t,\ \alpha \in \left (\alpha_{\ast},\alpha^{\ast}\right ],\ t\in
\lbrack 0,\left (\alpha -\alpha_{\ast}\right )b).
\end{equation*}
Obviously, $p_{b}(\alpha ,t)$ is decreasing in $t$ and increasing in $\alpha
$, and satisfies inequality $0<p_{b}(\alpha ,t)\leq 1$.

We introduce the space $M_{b}$ of square-integarble progessively measurable
random processes $u:\left
[0,\left
(\alpha^{\ast}-\alpha_{\ast}\right
)b%
\right )\rightarrow X_{\alpha^{\ast}}$ such that $u(t)\in X_{\alpha }$ for $%
t<\left (\alpha -\alpha_{\ast}\right )b$, and
\begin{equation*}
\left \vert \left \vert \left \vert u\right \vert \right \vert \right \vert
_{b}:=\sup \left \{\left ({\mathbb{E}}\left \Vert u(t)\right \Vert _{\alpha
}^{2}p_{b}(\alpha ,t)\right )^{1/2}:\ \alpha \in \left
(\alpha_{\ast},\alpha^{\ast}\right ],t\in \lbrack 0,\left (\alpha
-\alpha_{\ast}\right )b)\right \}<\infty .
\end{equation*}
Thus for any $u\in M_{b}$ there exists $C>0$ such that
\begin{equation*}
{\mathbb{E}}\left \Vert u(t)\right \Vert _{\alpha }^{2}\leq \frac{C}{1-\left
(\left (\alpha -\alpha_{\ast}\right )b\right )^{-1}t},\ t<\left (\alpha
-\alpha_{\ast}\right )b.
\end{equation*}
The pair $M_{b}$, $\left \vert \left \vert \left \vert \cdot \right \vert
\right \vert \right \vert _{b}$ forms a separable Banach space. For any $a>b$ there is a natural map $O_{ab}:M_{a}\rightarrow M_{b}$
given by the restriction
\begin{equation*}
O_{ab}u=u\upharpoonright _{\left[ 0,\left(
\alpha^{\ast}-\alpha_{\ast}\right) b\right) }.
\end{equation*}%
\begin{remark}
Similar spaces of deterministic functions $u:\left [0,\left (\alpha^{\ast}
-\alpha_{\ast}\right )b\right )\rightarrow X_{\alpha^{\ast}}$ where used in
\cite{Ni,Saf,BHP}.
\end{remark}

\begin{remark}
\label{Rem-restr}For any fixed $b>0,$ $T<\left
(\alpha^{\ast}-\alpha_{\ast}%
\right )b$ and $\beta \in \left (Tb^{-1}+ \alpha_{\ast},\alpha^{\ast}%
\right
]$ consider the spaces $E_{\beta ,T}$ and $H_{\beta ,T}$ of
square-integarble progessively measurable random processes $%
u:\left
[0,T\right )\rightarrow X_{\beta }$ and $h:\left
[0,T\right
)%
\rightarrow H_{\beta }$ with finite norms
\begin{equation*}
\left \Vert u\right \Vert _{E_{\beta ,T}}:=\sup _{t\in \left [0,T\right
)}\left ({\mathbb{E}}\left \Vert u(t)\right \Vert _{\beta }^{2}\right )^{1/2}%
{\text { and }}\left \Vert h\right \Vert _{H_{\beta ,T}}:=\sup _{t\in \left
[0,T\right )}\left ({\mathbb{E}}\left \Vert u(t)\right \Vert _{H_{\beta
}}^{2}\right )^{1/2},
\end{equation*}
respectively. Let $u^{(T)}:=u\upharpoonright _{\left [0,T\right )}$ be the
restriction of a process $u\in M_{b}$ to time interval $\left [0,T\right )$.
Observe that $p_{b}(\beta ,t)\geq c$ for some constant $c>0$ and all 
$t\leq T$. Thus $\left \Vert u^{T}\right \Vert
_{E_{\beta ,T}}\leq c^{-1}\left \vert \left \vert \left \vert u\right \vert
\right \vert \right \vert _{b}^{2}$ and so $u^{(T)}\in E_{\beta ,T}$.
Moreover, it is clear that for any $f\in {\mathcal{GL}}^{(1)}$ and $B\in {%
\mathcal{GL}}^{(2)}$ we have $f(u^{(T)})\in E_{\beta ,T}$ and $B(u^{(T)})\in
H_{\beta ,T}$. Indeed, we can fix $\alpha \in \left (Tb^{-1}+
\alpha_{\ast},\beta \right )$ (so that $u^{(T)}\in E_{\alpha ,T}$) and apply
estimates from Remark \ref{rem-bound}, which will show that 
$\left\Vert f(u))\right\Vert _{E_{\beta ,T}},\left\Vert B(u(t))\right\Vert _{H_{\beta,T }}<\infty$.
\end{remark}

From now on, we fix $f\in {\mathcal{GL}}^{(1)}$ and $B\in {\mathcal{GL}}%
^{(2)}$. For any $u\in M_{b}$ define
\begin{equation}
F(u)(t)=\int_{0}^{t}f(u(s))ds+\int_{0}^{t}B(u(s))dW(s)\in X_{\alpha ^{\ast
}},\ t<b\alpha ^{\ast }.  \label{lip-map}
\end{equation}%
According to Remark \ref{Rem-restr}, $f(u^{(t)}(\cdot ))\in E_{\beta ,t}$
and $B(u^{(t)}(\cdot ))\in HS_{\beta ,t}$ for any $\beta >b^{-1}t+\alpha
_{\ast }$. Thus the right-hand side is well-defined in $X_{\beta }$ with $%
\beta >b^{-1}t+\alpha _{\ast }$.

Consider equation
\begin{equation}
u=u_{0}+F(u)  \label{SDE}
\end{equation}%
with $u_{0}\in X_{\alpha ^{\ast }}$, cf. (\ref{SDE-int}), and set $b^{\ast
}:=\frac{\sqrt{1+(\alpha ^{\ast }-\alpha _{\ast })L^{-1}}-1}{2(\alpha ^{\ast
}-\alpha _{\ast })}$. The following theorem states the main existence result
of this paper.

\begin{theorem}[Existence]
\label{theor-main}Equation (\ref{SDE}) has a solution $u\in M_{b}$ for any $%
b<b^{\ast }$. It is unique in the following sense: if $u_{1}\in M_{b_{1}}$
and $u_{2}\in M_{b_{2}}$ are two solutions and $b_{1}\leq b_{2}<b^{\ast }$
than $O_{b_{2}b_{1}}u_{2}=u_{1}$.
\end{theorem}

\noindent \textbf{Proof.} It is sufficient to show that the map
\begin{equation*}
u\mapsto u_{0}+F(u)
\end{equation*}
is contractive in $M_{b}$ with $b<b^{\ast }$, which in turn will imply the
existence of its (unique) fixed point. It is straightforward that if $u\ $is
the fixed point in $M_{b_{1}}$ then $O_{b_{2}b_{1}}u$ is the fixed point in $%
M_{b_{2}}$. Thus the statement of the theorem follows from Theorem \ref%
{Lip-gen} and Corollary \ref{corr-contr}, which will be proved in the next
section. \hfill $\square $

Of course the choice of the weight function $p_{b}$ is somehow ambiguous.
The following statement is a corollary of Theorem \ref{theor-main}
formulated in a slightly more invariant form (although with some loss of
information).

\begin{corollary}
\label{corr-exu}Equation (\ref{SDE}) has a solution $u:\left[ 0,\left( \bar{%
\alpha}-\alpha_{\ast}\right) b^{\ast }\right) \rightarrow X_{\bar{\alpha}}$.
Moreover, $u\upharpoonright _{\left[ 0,T\right) }\in E_{\beta ,T}$ for any $%
T<\left( \alpha^{\ast}-\alpha_{\ast}\right) b^{\ast }$ and $\beta \in \left(
T/b^{\ast }+\alpha_{\ast},\alpha^{\ast}\right] $.
\end{corollary}

\bigskip

\noindent Theorem \ref{theor-main} establishes the uniqueness of the
solution in $M_{b}$. A natural question that arises here is whether there
might be a solution that does not belong to any $M_{b}$. An answer is given
by the following (somewhat stronger) uniqueness result.

\begin{theorem}[Uniqueness]
\label{theor-uniq}Fix $\beta \in \left[ \alpha_{\ast},\frac{\alpha_{\ast}
+\alpha^{\ast}}{2}\right] $ and $b<b^{\ast }$ and assume that $u\in E_{\beta
,T}$, where $T=\left( \alpha^{\ast}-\alpha_{\ast}\right) b$, is a solution
of equation (\ref{SDE}). Then $u\in M_{b}$ and coincides in this space with
the solution from Theorem \ref{theor-main}.
\end{theorem}

\noindent \textbf{Proof.} First observe that $E_{\alpha_{\ast},T}\subset
M_{b}$, which implies the statement for $\beta =\alpha_{\ast}$.

Let now $\beta \in (\alpha_{\ast},\alpha^{\ast})$ and us consider the Banach
space $M_{b,\beta }$ defined by replacing $\alpha_{\ast}$ with $\beta $ in
the definition of $M_{b}$ (so that $M_{b}=M_{b,\alpha_{\ast}}$). Then we
clearly have $OE_{\beta ,T}\subset M_{b,\beta }$, with the operator $O$
given by the restriction to time interval $\left[ 0,\left( \bar{\alpha}%
-\beta \right) b\right) $. Moreover, $OM_{b}\subset M_{b,\beta }$. Indeed,
for any $v\in M_{b}$ and $t\in \left[ 0,\left( \alpha^{\ast}-\beta \right)
b\right) $ we have $v(t)\in \bigcap\limits_{\alpha
>tb^{-1}+\alpha_{\ast}}X_{a}\subset \bigcap\limits_{\alpha >tb^{-1}+\beta
}X_{a} $ because $\beta >\alpha_{\ast}$. A direct check shows that $%
\left\vert \left\vert \left\vert u\right\vert \right\vert \right\vert
_{b}\geq \left\vert \left\vert \left\vert u\right\vert \right\vert
\right\vert _{b,\beta }$.

Observe that the proof of Theorem \ref{Lip-gen} (and thus of Theorem \ref%
{theor-main}) can be accomplished in the space $M_{b,\beta }$ instead of $%
M_{b}$, which implies that $Ou$ is the unique solution of (\ref{SDE}) in $%
M_{b,\beta }$. Let now $v\in M_{b}$ be the solution constructed in Theorem %
\ref{theor-main}. By the uniqueness part of that theorem, we have $Ou=Ov$,
which means that $u(t)=v(t)$, $t\in \left[ 0,\left( \alpha^{\ast}-\beta
\right) b\right) $. Observe that the assumption $\beta \leq \frac{%
\alpha_{\ast}+\alpha^{\ast}}{2}$ implies that $\left( \alpha^{\ast}-\beta
\right) b\geq \left( \beta -\alpha_{\ast}\right) b$. By Lemma \ref{lemma-un}
below we have $u\in M_{b}$, and the statement of the theorem follows from
the uniqueness in $M_{b}$. \hfill $\square $

\begin{lemma}
\label{lemma-un}Let $\beta \in (\alpha_{\ast},\alpha^{\ast})$, $u\in
E_{\beta ,T}$ and there exist $v\in M_{b}$ such that
\begin{equation}
u\upharpoonright _{\left[ 0,\left( \beta -\alpha_{\ast}\right) b\right)
}=v\upharpoonright _{\left[ 0,\left( \beta -\alpha_{\ast}\right) b\right) }.
\label{uvequ}
\end{equation}%
Then $u\in M_{b}$.
\end{lemma}

\noindent \textbf{Proof.} $u\in M_{b}$ iff $\exists C>0$ such that $\forall
\alpha \in (\alpha_{\ast},\alpha^{\ast})$ we have $u(t)\in X_{\alpha }$ for $%
t<\left( \alpha -\alpha_{\ast}\right) b$ and $\sup_{t<\left( \alpha
-\alpha_{\ast}\right) b}{\mathbb{E}}\left\Vert u(t)\right\Vert _{\alpha
}^{2}p_{b}(\alpha ,t)<C$. In our case, this holds for $\alpha <\beta $
because of (\ref{uvequ}) and for $\alpha \geq \beta $ because of the
inclusion $u\in E_{\beta ,T}$ and the bound $p_{b}(\alpha ,t)<1$.\hfill $%
\square $ \bigskip

\noindent Our main example is given by an infinite system of SDEs describing
stochastic dynamics of certain infinite particle spin system and will be
discussed in Section \ref{spins}. Here, we provide an example of a very
different type, which can also be dealt with by much simpler methods and
thus clarifies up to some extend the statement of Theorem \ref{theor-main}.

\begin{remark}\label{Banach}For simplicity, we required $X_{\alpha }$ to be Hilbert spaces. This is in
fact not essential and the case of a scale of suitable Banach spaces can be
treated in a similar way.
\end{remark}

\begin{example}
Consider the following SPDE on the $1$-dimensional torus ${\mathbb{T}}$:
\begin{equation}
du(t)=cu_{x}(t)dW(t),  \label{SPDE}
\end{equation}
where $u(t)\in C^{1}({\mathbb{T}})$, $u_{x}(t,x):=\frac{{\mathcal{\partial }}%
}{\partial x}u(t,x)$, $x\in {\mathbb{T}}$, $c\in {\mathbb{R}}$ and $W$ is a
real-valued Wiener process. Denote by $\widehat {v}(k)$, $k\in {\mathbb{Z}}$%
, the Fourier coefficients of $v\in L^{2}({\mathbb{T}})$ and define the
scale of Hilbert spaces
\begin{equation*}
X_{\alpha }:=\left \{v\in L^{2}({\mathbb{T}}):\left \Vert v\right \Vert
_{\alpha }:=\left (\sum _{k\in {\mathbb{Z}}}\left \vert \widehat {v}%
(k)\right \vert ^{2}e^{\alpha \left \vert k\right \vert ^{2}}\right
)^{1/2}<\infty \right \},\ \alpha >0.
\end{equation*}
It is clear that $X_{\alpha }\subset X_{\beta },\ \alpha >\beta $ (cf.
Remark \ref{rem-scale}). Let ${\mathcal{H}}:={\mathbb{R}}$ and define $%
B:X_{\alpha }\rightarrow HS({\mathcal{H}},X_{\beta })$ by the formula $%
B(v)h=bv_{x}h$, $v\in X_{\alpha },$ $h\in {\mathcal{H}}$. Equation (\ref%
{SPDE}) can now be written in the form (\ref{SDE-int}). Moreover, it can be
shown by a direct computation that $B$ satisfies condition (\ref{lg-diff}).
Thus, by Theorem \ref{theor-main} adopted to this setting, for any $\beta
<\alpha $ and an initial condition $u(0)\in X_{\alpha }$ there exists a
solution $u(t)\in X_{\beta }$, $t<\tau (\alpha -\beta )$, where $\tau $ is a
constant (independent of $\alpha $ and $\beta $ but possibly dependent on
their allowed range).

Observe that equation (\ref{SPDE}) can be solved explicitly. Indeed, the
Fourier coefficients of $u(t)$ satisfy the equation
\begin{equation*}
d\widehat {u}(t,k)=ick\widehat {u}(t,k)dW(t),\ k\in {\mathbb{Z}},
\end{equation*}
so that
\begin{equation*}
\widehat {u}(t,k)=e^{tc^{2}k^{2}/2}e^{ickW(t)}\widehat {u}(0,k),\ k\in {%
\mathbb{Z}},
\end{equation*}
which in turn implies the equality
\begin{equation}
\left \vert \widehat {u}(t,k)\right \vert ^{2}=e^{tc^{2}k^{2}}\left \vert
\widehat {u}(0,k)\right \vert ^{2},\ k\in {\mathbb{Z}}.  \label{norm}
\end{equation}
Fix any $\beta <\alpha $ and an initial condition $u(0)\in X_{\alpha }$. It
follows directly from (\ref{norm}) that the solution $u(t)\ $belongs to $%
X_{\beta }$ for $t<c^{-2}(\alpha -\beta )$. It is also clear that the
solution does not live in the scale of standard Sobolev spaces. Neither of
course does $B$ satisfy condition (\ref{lg-diff}) in such a scale.
\end{example}

\section{Proof of the contractivity.}

In this section, we will show that $F$ is a contraction in $M_{b}$ with $b$
sufficiently small.

\begin{theorem}
\label{Lip-gen}For any $b>0$, formula (\ref{lip-map}) defines the map $%
F:M_{b}\rightarrow M_{b}$. Moreover, $F$ is Lipschitz continuous with
Lipschitz constant $2bL\sqrt{\left( \alpha^{\ast}-\alpha_{\ast}\right)
+b^{-1}}$.
\end{theorem}

\noindent \textbf{Proof. }Let $u,v\in M_{b}$ and fix $\beta \leq
\alpha^{\ast} $ and $t\in (0,b\beta )$. Then $F(u)(t),F(v)(t)\in X_{\beta }$%
, and we have the estimate
\begin{multline*}
{\mathbb{E}}\left\Vert F(u)(t)-F(v)(t)\right\Vert _{\beta }^{2}\leq t{%
\mathbb{E}}\int_{0}^{t}\left\Vert f(u(s))-f(v(s))\right\Vert _{\beta }^{2}ds
\\
+{\mathbb{E}}\int_{0}^{t}\left\Vert B(u(s))-B(v(s))\right\Vert _{H_{\beta
}}^{2}ds \\
\leq cL^{2}{\mathbb{E}}\int_{0}^{t}\left\Vert u(s)-v(s)\right\Vert _{\alpha
(s)}^{2}\left( \beta -\alpha (s)\right) ^{-1}ds
\end{multline*}%
with $c=(b\left( \alpha^{\ast}-\alpha_{\ast}\right) +1)$, for any $\alpha
(s) $ satisfying $b^{-1}s+\alpha_{\ast}<\alpha (s)<\beta $. Then
\begin{multline}
{\mathbb{E}}\left\Vert F(u)(t)-F(v)(t)\right\Vert _{\beta }^{2}\leq cL^{2}{%
\mathbb{E}}\int_{0}^{t}\left\Vert u(s)-v(s)\right\Vert _{\alpha
(s)}^{2}p_{b}(\alpha (s),s)  \label{i1} \\
\times p_{b}(\alpha (s),s)^{-1}\left( \beta -\alpha (s)\right) ^{-1}ds \\
\leq cL^{2}\left\vert \left\vert \left\vert u-v\right\vert \right\vert
\right\vert _{b}^{2}\int_{0}^{t}p_{b}(\alpha (s),s)^{-1}\left( \beta -\alpha
(s)\right) ^{-1}ds.
\end{multline}%
We set
\begin{equation*}
\alpha (s)=\frac{1}{2}\left( \beta +b^{-1}s+\alpha_{\ast}\right) .
\end{equation*}%
Then
\begin{equation*}
\beta -\alpha (s)=\frac{1}{2}\left( \hat{\beta}-b^{-1}s\right) ,\ \hat{\beta}%
:=\beta -\alpha_{\ast},
\end{equation*}%
and
\begin{equation*}
p_{b}(\alpha (s),s)=\left( \hat{\beta}-b^{-1}s\right) \left( \hat{\beta}%
+b^{-1}s\right) ^{-1},
\end{equation*}%
and the integral term of (\ref{i1}) obtains the form
\begin{multline*}
I:=2\int_{0}^{t}\left( \hat{\beta}-b^{-1}s\right) ^{-2}\left( \hat{\beta}%
+b^{-1}s\right) ds \\
\leq 2b\left[ \left( \hat{\beta}-b^{-1}t\right) ^{-1}-\hat{\beta}^{-1}\right]
\left( \hat{\beta}+b^{-1}t\right) \\
\leq 2b\left( \hat{\beta}-b^{-1}t\right) ^{-1}\hat{\beta}\left( 1+\hat{\beta}%
^{-1}b^{-1}t\right) \\
=2bp_{b}(\beta ,t)^{-1}\left( 1+\hat{\beta}^{-1}b^{-1}t\right) .
\end{multline*}%
The bound $\hat{\beta}^{-1}b^{-1}t<1$ implies that
\begin{equation*}
I\leq 4bp_{b}(\beta ,t)^{-1}.
\end{equation*}%
Thus it follows from (\ref{i1}) that
\begin{equation}
\left\vert \left\vert \left\vert F(u)-F(v)\right\vert \right\vert
\right\vert _{b}\leq 2\sqrt{c}L\left\vert \left\vert \left\vert
u-v\right\vert \right\vert \right\vert _{b}.  \label{Lip1}
\end{equation}

Let us now show that $F$ preserves the space $M_{b}$. For this, we set ${%
\mathbf{u}}_{0}(t)=0\in X_{\alpha_{\ast}}$. Then ${\mathbf{u}}_{0}\in M_{b}$
so that $F(u)-F({\mathbf{u}}_{0})\in M_{b}$ provided $u\in M_{b}$. Moreover,
\begin{equation*}
F({\mathbf{u}}_{0})(t)=tf(0)+B(0)W(t),
\end{equation*}%
and so
\begin{equation*}
{\mathbb{E}}\left\Vert F({\mathbf{u}}_{0})(t)\right\Vert _{\beta }^{2}\leq
2t^{2}\left\Vert f(0)\right\Vert _{\beta }^{2}+2t\left\Vert B(0)\right\Vert
_{H_{\beta }}^{2}\leq 2(t^{2}+t)K^{2}\beta ^{-1}.
\end{equation*}%
In the second inequality we used Remark \ref{rem-bound} with $u=0$ and $%
\alpha =0$. Then
\begin{equation*}
\left\vert \left\vert \left\vert F({\mathbf{u}}_{0})\right\vert \right\vert
\right\vert _{b}^{2}\leq \sup_{\beta ,t:~t<b(\beta
-\alpha_{\ast})}p_{b}(\beta ,t)2(t^{2}+t)K^{2}\beta ^{-1}\leq 2cK^{2}<\infty
,
\end{equation*}%
because $p_{b}(\beta ,t)\leq 1$ and $t<b\beta \leq b\alpha^{\ast}$. Thus $F({%
\mathbf{u}}_{0})\in M_{b}$ and
\begin{equation*}
F(u)=\left( F(u)-F(u_{0})\right) +F(u_{0})\in M_{b}.
\end{equation*}%
This together with (\ref{Lip1}) implies the result. \hfill $\square $

\begin{corollary}
\label{corr-contr}The map $F$ is contractive in every $M_{b}$ with $b<\frac{%
\sqrt{1+(\alpha^{\ast}-\alpha_{\ast})L^{-1}}-1}{2(\alpha^{\ast}-\alpha_{%
\ast})}$.
\end{corollary}

\section{Stochastic spin dynamics of a quenched particle system \label{spins}%
}

Our main example is motivated by the study of stochastic dynamics of
interacting particle systems. Let $\gamma \subset X={\mathbb{R}}^{d}$ be a
locally finite set (configuration) representing a collection of point
particles. Each particle with position $x\in X$ is characterized by an
internal parameter (spin) $\sigma _{x}\in S={\mathbb{R}}^{1}$.

We fix a configuration $\gamma $ and look at the time evolution of spins $%
\sigma _{x}(t)$, $x\in \gamma $, which is described by a system of
stochastic differential equations in $S$ of the form
\begin{equation}
d\sigma _{x}(t)=f_{x}(\bar{\sigma})dt+B_{x}(\bar{\sigma})dW_{x}(t),\ x\in
\gamma ,  \label{system}
\end{equation}%
where $\bar{\sigma}=(\sigma _{x})_{x\in \gamma }$ and $W=(W_{x})_{x\in
\gamma }$ is a collection of independent Wiener processes in $S$. We assume
that both drift and diffusion coefficients $f_{x}$ and $B_{x}$ depend only
on spins $\sigma _{y}$ with $\left\vert y-x\right\vert <r$ for some fixed
interaction radius $r>0$ and have the form
\begin{equation}
f_{x}(\bar{\sigma})=\sum_{y\in \gamma }\varphi _{xy}(\sigma _{x},\sigma
_{y}),\ \ B_{x}(\bar{\sigma})=\sum_{y\in \gamma }\Psi _{xy}(\sigma
_{x},\sigma _{y}),  \label{system-coef}
\end{equation}%
where the mappings $\varphi _{xy}:S\times S\rightarrow S$ and $\Psi
_{xy}:S\times S\rightarrow S$ satisfy finite range and uniform Lipschitz
conditions, see Definition \ref{def-admiss} and Condition \ref{cond-inter1}
below.

Our aim is to realise system (\ref{system}) as an equation in a suitable
scale of Hilbert spaces and apply the results of previous sections in order
to find its strong solutions.

We introduce the following notations:

- $S^{\gamma }:=\prod _{x\in \gamma }S_{x}\ni \bar {\sigma }=(\sigma
_{x})_{x\in \gamma },\ \sigma _{x}\in S_{x}=S$;

- $\gamma _{x,r}:=\left \{y\in \gamma :\left \vert x-y\right \vert
<r\right
\},\ x\in \gamma $;

- $n_{x}\equiv n_{x,r}(\gamma ):=$ number of points in $\gamma _{x,r}$ ( $=$
number of particles interacting with particle in position $x$).

Observe that, although the number $n_{x}$ is finite, it is in general
unbounded function of $x$. We assume that it satisfies the following
regularity condition.

\begin{condition}
\label{cond-gamma}There exists a constant $a(\gamma ,r)$ such that
\begin{equation}
n_{x,r}(\gamma )\leq a(\gamma ,r)\left( 1+\left\vert x\right\vert \right)
^{1/2}  \label{c-gamma}
\end{equation}%
for all $x\in X$.
\end{condition}

\begin{remark}
Condition (\ref{c-gamma}) holds if $\gamma $ is a typical realization of a
Poisson or Gibbs (Ruelle) point process in $X$. For such configurations,
stronger (logarithmic) bound holds:
\begin{equation*}
n_{x,r}(\gamma )\leq c(\gamma )\left[ 1+\log (1+\left\vert x\right\vert )%
\right] r^{d},
\end{equation*}%
see e.g. \cite{R70} and \cite[p.~1047]{K93}.
\end{remark}

\subsection{Existence of the dynamics}

Our dynamics will live in the scale of Hilbert spaces
\begin{equation*}
X_{\alpha }=S_{\alpha }^{\gamma }:=\left\{ \bar{q}\in S^{\gamma }:\left\Vert
\bar{q}\right\Vert _{\alpha }:=\sqrt{\sum_{x\in \gamma }\left\vert
q_{x}\right\vert ^{2}e^{-\alpha \left\vert x\right\vert }}<\infty \right\}
,\ \alpha >0.
\end{equation*}%
Let us define the corresponding spaces ${\mathcal{GL}}^{(1)}$ and ${\mathcal{%
GL}}^{(2)}$ (cf. Condition \ref{Lip-drift}) and set
\begin{equation*}
{\mathcal{H}}=S_{0}^{\gamma }:=\left\{ \bar{q}\in S^{\gamma }:\left\Vert
\bar{q}\right\Vert _{0}:=\sqrt{\sum_{x\in \gamma }\left\vert
q_{x}\right\vert ^{2}}<\infty \right\} .
\end{equation*}%
Observe that $W(t):=\left( W_{x}(t)\right) _{x\in \gamma }$ is a cylinder
Wiener process in $\mathcal{H}$.

\bigskip

Let ${\mathcal{V}}$ be a family of mappings $V_{xy}:S^{2}\rightarrow S$, $%
x,y\in \gamma $.

\begin{definition}
\label{def-admiss}We call the family $\mathcal{V}$ admissible if it
satisfies the following two assumptions:
\end{definition}

\begin{itemize}
\item finite range: there exists constant $r>0$ such that $V_{xy}\equiv 0$
if $\left\vert x-y\right\vert \geq r$;

\item uniform Lipschitz continuity: there exists constant $C>0$ such that
\begin{equation}
\left\vert V_{xy}(q_{1}^{\prime },q_{2}^{\prime })-V_{xy}(q_{1}^{\prime
\prime },q_{2}^{\prime \prime })\right\vert \leq C\left( \left\vert
q_{1}^{\prime }-q_{1}^{\prime \prime }\right\vert +\left\vert q_{2}^{\prime
}-q_{2}^{\prime \prime }\right\vert \right)  \label{LipS2}
\end{equation}%
for all $x,y\in \gamma $ and $q_{1}^{\prime},q_{2}^{\prime},q_{1}^{\prime \prime },q_{2}^{\prime \prime }\in S$.
\end{itemize}

Define a map $\overline{V}:S^{\gamma }\rightarrow S^{\gamma }$ and a linear
operator $\widehat{V}(\bar{q}):S^{\gamma }\rightarrow S^{\gamma }$, $\bar{q}%
\in S^{\gamma }$, by the formula
\begin{equation*}
\overline{V}_{x}(\bar{q})=\sum_{y\in \gamma }V_{xy}(q_{x},q_{y}),
\end{equation*}%
and
\begin{equation*}
\left( \widehat{V}(\bar{q}){\bar{\sigma}}\right) _{x}:=\overline{V}_{x}(\bar{%
q})\sigma _{x},\,x\in \gamma ,\,{\bar{\sigma}}\in S^{\gamma },
\end{equation*}%
respectively.

\begin{lemma}
\label{lemma-gl}Assume that ${\mathcal{V}}$ is admissible. Then $\overline{V}%
\in {\mathcal{GL}}^{(1)}$ and $\widehat{V}\in {\mathcal{GL}}^{(2)}$.
\end{lemma}

The proof of this Lemma is quite tedious and will be given in Section \ref%
{sec-aux}.

\bigskip

Now we can return to the discussion of system (\ref{system}). Assume that
the following condition holds.

\begin{condition}
\label{cond-inter1} The families of mappings $\left\{ \varphi _{xy}\right\}
_{x,y\in \gamma }$ and $\left\{ \Psi _{xy}\right\} _{x,y\in \gamma }$ from (%
\ref{system-coef}) are admissible.
\end{condition}

By Lemma \ref{lemma-gl} we have $\overline{\varphi }\in {\mathcal{GL}}^{(1)}$
and $\widehat{\Psi }\in {\mathcal{GL}}^{(2)}$. Thus we can write (\ref%
{system}) in the form
\begin{equation*}
\bar{\sigma}(t)=\overline{\varphi }(\bar{\sigma})dt+\widehat{\Psi }(\bar{%
\sigma})dW(t),
\end{equation*}%
where $W(t)=\left( W_{x}(t)\right) _{x\in \gamma }$, and apply the results
of Section \ref{sec-main} to its integral counterpart. We summarize the
existence results in the following theorem, which follows directly from
Theorem \ref{theor-main}.

\begin{theorem}
\label{theor-particles}System (\ref{system}) has a strong solution $u:\left[
0,\left( \alpha ^{\ast }-\alpha _{\ast }\right) b^{\ast }\right) \rightarrow
X_{\alpha ^{\ast }}$. Moreover, $u(T)\in \bigcap\limits_{\alpha >T/b^{\ast
}+\alpha _{\ast }}
X_{\alpha }$ for any $T<\left( \alpha ^{\ast }-\alpha _{\ast }\right)
b^{\ast }$, and the restriction of $u$ to the time interval $\left[
0,T\right) $ belongs to $M_b$ with $b={\left( \alpha ^{\ast }-\alpha _{\ast }\right)^{-1}T}$.
\end{theorem}

\begin{remark} Theorem \ref{theor-particles} can also
be proved in the scale of Banach spaces $$S_{\alpha ,p}^{\gamma }:=\left\{ \bar{q}\in
S^{\gamma }:\left\Vert \bar{q}\right\Vert _{\alpha }:=\left( \sum_{x\in
\gamma }\left\vert q_{x}\right\vert ^{p}e^{-\alpha \left\vert x\right\vert
}\right) ^{1/p}<\infty \right\} ,\ \alpha >0,\ p>2,$$ cf. Remark \ref{Banach}. 
\end{remark}

\subsection{The uniqueness}

In this section we establish a stronger uniqueness result, extending to our
situation the method applied to deterministic systems in \cite{LLL}, \cite%
{DaF}. As before, the main ingredients here are the bound on the density of
configuration $\gamma $ (Condition \ref{cond-gamma}) and uniform Lipschitz
continuity of the maps $\varphi _{xy}$ and $\Psi _{xy}$ (Condition \ref%
{cond-inter1}). However, in contrast to the previous section, we will
consider solutions of a more general type.

\smallskip Let $E(S,T)$ be the space of square-integrable progressively
measurable random processes $q:[0,T)\rightarrow S$ such that $\sup _{t\in
\left [0,T\right
)}{\mathbb{E}}\left \Vert u(t)\right \Vert _{\beta
}^{2}<\infty$.

\begin{definition}
We call a random process $\bar{q}:[0,T)\rightarrow S^{\gamma }$ a pointwise
(strong) solution of system (\ref{system}) if $q_{x}(\cdot )\in E(S,T)$ and
satisfies integral equation
\begin{equation*}
q_{x}(t)=q_{x}(0)+\int_{0}^{t}f_{x}(\bar{q}(s))ds+\int_{0}^{t}B_{x}(\bar{q}%
(s))dW_{x}(s)
\end{equation*}%
for each $x\in \gamma $.
\end{definition}

It is clear that the solution constructed in Theorem \ref{theor-particles}
is a pointwise strong solution.

\begin{theorem}
\label{theor-point-uniq}Assume that Conditions \ref{cond-gamma} and \ref%
{cond-inter1} hold and let $\bar{q}^{(1)}(t),\bar{q}^{(2)}(t)\in S_{\beta
}^{\gamma }$ be two pointwise strong solutions of (\ref{system}) on $\left[
0,T\right) $, and let $\bar{q}^{(1)}(0)=\bar{q}^{(2)}(0)$ a.s. Then $\bar{q}%
^{(1)}(t)=\bar{q}^{(2)}(t)$ a.s. for any $t\in \lbrack 0,T)$.
\end{theorem}

To proceed with the proof, we need the following Lemma, which will in turn
be proved in Section \ref{sec-aux}. For any $n\in \mathbb{N}$ and$\ t\in
\lbrack 0,T)$ define%
\begin{equation*}
\delta _{n}(t):=\sup_{\left\vert x\right\vert \leq nr}\mathbb{E}\left\vert
q_{x}^{(1)}(t)-q_{x}^{(2)}(t)\right\vert ^{2}.
\end{equation*}

\begin{lemma}
\label{lem-unique}Assume that conditions of Theorem \ref{theor-point-uniq}
hold. Then there exists $\mu >0$ such that%
\begin{equation}
\mathbb{\delta }_{n}(t)\leq 2n(t+1)\mu \int_{0}^{t}\delta _{n+1}(s)ds
\label{delta-iter}
\end{equation}%
for any $t\in \lbrack 0,T)$.
\end{lemma}

\noindent \textbf{Proof of Theorem \ref{theor-point-uniq}.} The $N$-th
iteration of bound (\ref{delta-iter}) gives the estimate
\begin{equation}
\mathbb{\delta }_{n}(t)\leq \frac{\left( 2(t+1)t\mu \right) ^{N}}{N!}%
n(n+1)....(n+N-1)\sup_{s\leq t}\delta _{n+N}(s)  \label{estim2}
\end{equation}%
for any $N=2,3,...$ . Set
\begin{equation*}
R:=\sup_{s\leq T}\left\{ \mathbb{E}\left\Vert \bar{q}^{(1)}(s)\right\Vert
_{\beta }^{2},\mathbb{E}\left\Vert \bar{q}^{(2)}(s)\right\Vert _{\beta
}^{2}\right\} .
\end{equation*}%
Taking into account that $\bar{q}^{(1)}(t),\bar{q}^{(2)}(t)\in S_{\beta
}^{\gamma }$ we obtain the bounds
\begin{equation*}
\mathbb{E}\left\vert q_{x}^{(i)}(t)\right\vert ^{2}\leq e^{\beta \left\vert
x\right\vert }\mathbb{E}\left\Vert \bar{q}^{(i)}(t)\right\Vert _{\beta
}^{2}\leq e^{\beta \left\vert x\right\vert }R,\,i=1,2,
\end{equation*}%
which imply that
\begin{equation*}
\delta _{n+N}(s)\leq 4e^{\beta (n+N)r}R
\end{equation*}%
for any $s\in \lbrack 0,T]$. It follows now from (\ref{estim2}) that%
\begin{multline*}
\mathbb{\delta }_{n}(t)\leq 4e^{\beta (n+N)r}R\frac{\left( 2(t+1)t\mu
\right) ^{N}}{N!}n(n+1)....(n+N-1) \\
=4e^{\beta (n+N)r}R\left( 2(t+1)t\mu \right) ^{N}\left(
\begin{array}{c}
n+N-1 \\
N%
\end{array}%
\right)  \\
=4e^{\beta nr}R\left[ \left( 2e^{\beta r+1}\mu (t+1)t\right) \frac{n+N-1}{N}%
\right] ^{N}.
\end{multline*}%
Here we used the well-known inequality $\binom{M}{N}\leq \bigl(\frac{M\,e}{N}%
\bigr)^{N}$, $1\leq N\leq M$. For $N>n-1$ we have $\frac{n+N-1}{N}<2$ and so%
\begin{equation*}
\mathbb{\delta }_{n}(t)<4e^{\beta nr}R\left[ 4e^{\beta r+1}\mu (t+1)t\right]
^{N}\rightarrow 0,\ N\rightarrow \infty ,
\end{equation*}%
provided $4e^{\beta r+1}\mu (t+1)t<1$ (e.g. $t<t_{0}:=\frac{1}{4}\left(
e^{\beta r+1}\mu (\alpha ^{\ast }+1)b\right) ^{-1}$). Thus%
\begin{equation*}
\sup_{\left\vert x\right\vert \leq nr}\mathbb{E}\left\vert
q_{x}^{(1)}(t)-q_{x}^{(2)}(t)\right\vert ^{2}=0,\ t<t_{0},
\end{equation*}%
for all $n\geq 1$, so that $\bar{q}^{(1)}(t)=\bar{q}^{(2)}(t)$ a.s. for any $%
t\in \lbrack 0,t_{0})$.

These arguments can be repeated on each of the time intervals $\left[
t_{k},t_{k+1}\right) $ with $t_{k}:=kt_{0},\ k=1,2,...,$ which shows that $%
\bar{q}^{(1)}(t)=\bar{q}^{(2)}(t)$ a.s. for any $t\in \lbrack 0,T)$, and the
proof is complete.
\hfill%
$\square $

\section{Proofs of auxiliary results}\label{sec-aux}

In this section, we present proofs of two technical lemmas used in the previous section.

\subsection{Proof of Lemma \protect\ref{lemma-gl}}

\textbf{Step 1. }We first show that $\overline{V}$ is a mapping $S_{\alpha
}^{\gamma }\rightarrow S_{\beta }^{\gamma }$ for any $\alpha <\beta $. For
any $\bar{q}\in S_{\alpha }^{\gamma }$ we have%
\begin{eqnarray*}
\left\Vert \overline{V}(\bar{q})\right\Vert _{\beta }^{2} &=&\sum_{x\in
\gamma }\left\vert \sum_{y\in \gamma }V_{xy}(q_{x},q_{y})\right\vert
^{2}e^{-\beta \left\vert x\right\vert } \\
&\leq &3C^{2}\sum_{x\in \gamma }\sum_{y\in \gamma _{x,r}}n_{x}\left(
1+\left\vert q_{x}\right\vert ^{2}+\left\vert q_{y}\right\vert ^{2}\right)
e^{-\beta \left\vert x\right\vert }.
\end{eqnarray*}%
The polynomial bound on the growth of $n_{x}$ implies that
\begin{equation*}
\sum_{x\in \gamma }\sum_{y\in \gamma _{x,r}}n_{x}e^{-\beta \left\vert
x\right\vert }=\sum_{x\in \gamma }n_{x}^{2}e^{-\beta \left\vert x\right\vert
}\leq \sum_{x\in \gamma }n_{x}^{2}e^{-\alpha _{\ast }\left\vert x\right\vert
}=:c(\gamma ,\alpha _{\ast })<\infty .
\end{equation*}%
Next, we estimate%
\begin{multline*}
\sum_{x\in \gamma }\sum_{y\in \gamma _{x,r}}n_{x}\left\vert q_{x}\right\vert
^{2}e^{-\beta \left\vert x\right\vert }=\sum_{x\in \gamma
}n_{x}^{2}\left\vert q_{x}\right\vert ^{2}e^{-\left( \beta -\alpha \right)
\left\vert x\right\vert }e^{-\alpha \left\vert x\right\vert } \\
\leq \sup_{x\in \gamma }\left( n_{x}^{2}e^{-\left( \beta -\alpha \right)
\left\vert x\right\vert }\right) \left\Vert \bar{q}\right\Vert _{\alpha
}^{2}.
\end{multline*}%
Observe that $\sum\limits_{x\in \gamma }\sum\limits_{y\in \gamma
_{x,r}}=\sum\limits_{\substack{ x,y\in \gamma  \\ \left\vert x-y\right\vert
<r}}=\sum\limits_{y\in \gamma }\sum\limits_{x\in \gamma _{y,r}}$, and so
\begin{multline*}
\sum_{x\in \gamma }\sum_{y\in \gamma _{x,r}}n_{x}\left\vert q_{y}\right\vert
^{2}e^{-\beta \left\vert x\right\vert }\leq e^{\beta r}\sum_{y\in \gamma
}N_{y}\left\vert q_{y}\right\vert ^{2}e^{-(\beta -\alpha )\left\vert
y\right\vert }e^{-\alpha \left\vert y\right\vert } \\
\leq e^{\beta r}\sup_{y\in \gamma }\left( N_{y}e^{-(\beta -\alpha
)\left\vert y\right\vert }\right) \left\Vert \bar{q}\right\Vert _{\alpha
}^{2},
\end{multline*}%
where $N_{y}:=\sum_{x\in \gamma _{y,r}}n_{x}$. Here we used inequality $%
\left\vert y\right\vert \leq \left\vert y-x\right\vert +\left\vert
x\right\vert \leq r+\left\vert x\right\vert $ for $y\in \gamma _{x,r}$, so
that $e^{-\beta \left\vert x\right\vert }\leq e^{\beta r}e^{-\beta
\left\vert y\right\vert }$. Condition \ref{cond-gamma} implies that%
\begin{equation*}
N_{x}\leq a(\gamma ,r)^{2}\left( 1+\left\vert x\right\vert \right)
^{1/2}\left( 1+r+\left\vert x\right\vert \right) ^{1/2}<a(\gamma
,r)^{2}(1+r)^{1/2}\left( 1+\left\vert x\right\vert \right) ,
\end{equation*}%
and
\begin{equation*}
n_{x}^{2}\leq a(\gamma ,r)^{2}\left( 1+\left\vert x\right\vert \right)
\end{equation*}%
for any $x\in \gamma $. Setting $c_{2}(\gamma ,r):=a(\gamma ,r)^{2}\left[
1+e^{\alpha ^{\ast }r}(1+r)^{1/2}\right] $ and $L^{2}=3C^{2}\left(
c_{1}+c_{2}\right) e^{\alpha ^{\ast }-\alpha _{\ast }-1}$ we obtain the
bound
\begin{equation*}
\left\Vert \overline{V}(\bar{q})\right\Vert _{\beta }^{2}\leq 3C^{2}\left(
c_{1}+c_{2}\right) \left[ \sup_{s>0}(1+s)e^{-(\beta -\alpha )s}\right]
\left\Vert \bar{q}\right\Vert _{\alpha }^{2}\leq L^{2}\left( \beta -\alpha
\right) ^{-1}\left\Vert \bar{q}\right\Vert _{\alpha }^{2}<\infty .
\end{equation*}

\textbf{Step 2. }Lipschitz condition (\ref{LipS2}) implies the estimate
\begin{eqnarray*}
\left\Vert \overline{V}(\bar{q}^{\prime })-\overline{V}(\bar{q}^{\prime
\prime })\right\Vert _{\beta }^{2} &=&\sum_{x\in \gamma }\left\vert
\sum_{y\in \gamma }V_{xy}(q_{x}^{\prime },q_{y}^{\prime })-\sum_{y\in \gamma
}V_{xy}(q_{x}^{\prime \prime },q_{y}^{\prime \prime })\right\vert
^{2}e^{-\beta \left\vert x\right\vert } \\
&\leq &2C^{2}\sum_{x\in \gamma }\sum_{y\in \gamma _{x,r}}n_{x}\left(
\left\vert q_{x}^{\prime }-q_{x}^{\prime \prime }\right\vert ^{2}+\left\vert
q_{y}^{\prime }-q_{y}^{\prime \prime }\right\vert ^{2}\right) e^{-\beta
\left\vert x\right\vert }
\end{eqnarray*}%
for any\textbf{\ }$\bar{q}^{\prime },\bar{q}^{\prime \prime }\in S_{\alpha
}^{\gamma }$. Similar to Step 1, we obtain the bound
\begin{multline*}
\left\Vert \overline{V}(\bar{q}^{\prime })-\overline{V}(\bar{q}^{\prime
\prime })\right\Vert _{\beta }^{2}\leq 2C^{2}c_{2}\left[
\sup_{s>0}(1+s)e^{-(\beta -\alpha )s}\right] \left\Vert \bar{q}^{\prime }-%
\bar{q}^{\prime \prime }\right\Vert _{\alpha }^{2} \\
\leq L^{2}\left( \beta -\alpha \right) ^{-1}\left\Vert \bar{q}^{\prime }-%
\bar{q}^{\prime \prime }\right\Vert _{\alpha }^{2}<\infty .
\end{multline*}

\textbf{Step 3.} 
The inclusion $\overline{V}(\bar{q})\in S_{\beta }^{\gamma }$ implies  that $\widehat{V}(\bar{q}){\bar{\sigma}}\in S_{\beta }^{\gamma }$
for any ${\bar{\sigma}}\in {\cal H}= S_{0}^{\gamma }$. A direct calculation shows that
$\widehat{V}(\bar{q}):{\cal H}\rightarrow {S}_{\beta
}^{\gamma }$ is a Hilbert-Schmidt operator with the norm equal to $%
\left\Vert \bar{V}(\bar{q})\right\Vert _{\beta }$. Thus the inclusion $%
\overline{V}\in {\mathcal{GL}}^{(1)}$ implies  that $\widehat{V}\in {%
\mathcal{GL}}^{(2)}$. \hfill $\square $

\subsection{Proof of Lemma \protect\ref{lem-unique}}

We start with the estimate of the distance between $q_{x}^{(1)}(t)$ and $%
q_{x}^{(2)}(t)$ for a fixed $x\in \gamma $ and $t\in \lbrack 0,T)$. From (%
\ref{system}) we obtain
\begin{multline}
\left\vert q_{x}^{(1)}(t)-q_{x}^{(2)}(t)\right\vert ^{2}\leq
2t\int_{0}^{t}\left\vert f_{x}(\bar{q}^{(1)}(s))-f_{x}(\bar{q}%
^{(2)}(s))\right\vert ^{2}ds  \label{estim11} \\
+2\int_{0}^{t}\left\vert B_{x}(\bar{q}^{(1)}(s))-B_{x}(\bar{q}%
^{(2)}(s))\right\vert^{2}ds=:2tI_{1,x}(t)+2I_{2,x}(t),
\end{multline}%
where $I_{1,x}(t)$ and $I_{2,x}(t)$ denote the first and second integral
terms, respectively. Taking into account that
\begin{equation*}
f_{x}(\bar{\sigma})=\sum_{y\in \gamma _{x}}\varphi _{xy}(\sigma _{x},\sigma
_{y}),\ \ B_{x}(\bar{\sigma})=\sum_{y\in \gamma }\Psi _{xy}(\sigma
_{x},\sigma _{y})
\end{equation*}%
and using Condition \ref{cond-inter1} we obtain
\begin{multline*}
I_{1,x}(t)\leq \int_{0}^{t}\left\vert \sum_{y\in \gamma _{x}}\left( \varphi
_{xy}(q_{x}^{(1)}(s),q_{y}^{(1)}(s))-\varphi
_{xy}(q_{x}^{(2)}(s),q_{y}^{(2)}(s))\right) \right\vert ^{2}ds \\
\leq n_{x}\int_{0}^{t}\sum_{y\in \gamma _{x}}\left\vert \varphi
_{xy}(q_{x}^{(1)}(s),q_{y}^{(1)}(s))-\varphi
_{xy}(q_{x}^{(2)}(s),q_{y}^{(2)}(s))\right\vert ^{2}ds \\
\leq 2n_{x}C^{2}\int_{0}^{t}\sum_{y\in \gamma _{x}}\left[ \left\vert
q_{x}^{(1)}(s)-q_{x}^{(2)}(s)\right\vert ^{2}+\left\vert
q_{y}^{(1)}(s)-q_{y}^{(2)}(s)\right\vert ^{2}\right] ds.
\end{multline*}%
Recall that
\begin{equation*}
n_{x}\leq a(\gamma ,r)\left( 1+\left\vert x\right\vert \right) ^{1/2}.
\end{equation*}%
Then for $\left\vert x\right\vert \leq nr$
\begin{multline*}
\mathbb{E}\left( I_{1,x}(t)\right) \leq 4n_{x}C^{2}\int_{0}^{t}\sum_{y\in
\gamma _{x}}\delta _{n+1}(s)ds=4n_{x}^{2}C^{2}\int_{0}^{t}\delta _{n+1}(s)ds
\\
\leq 4C^{2}a(\gamma ,r)^{2}\left( 1+\left\vert x\right\vert \right)
\int_{0}^{t}\delta _{n+1}(s)ds\leq \mu n\int_{0}^{t}\delta _{n+1}(s)ds
\end{multline*}%
with $\mu :=4C^{2}a(\gamma ,r)^{2}\left( 1+r\right) $. Similarly,
\begin{equation*}
\mathbb{E}\left( I_{2,x}(t)\right) \leq \mu n\int_{0}^{t}\delta _{n+1}(s)ds,
\end{equation*}%
so that (\ref{estim11}) implies the inequality
\begin{equation*}
\mathbb{E}\left\vert q_{x}^{(1)}(t)-q_{x}^{(2)}(t)\right\vert ^{2}\leq
2(t+1)\mu n\int_{0}^{t}\delta _{n+1}(s)ds
\end{equation*}%
and, consequently,
\begin{equation*}
\mathbb{\delta }_{n}(t)\leq 2(t+1)\mu n\int_{0}^{t}\delta _{n+1}(s)ds.
\end{equation*}%
The proof is complete.
\hfill%
$\square $

\bigskip

\textbf{Acknowledgment} \bigskip

I am deeply indebted to Yuri Kondratiev for his influence and kind support.
Stimulating discussions with Zdzislaw Brzezniak, Dmitri Finkelshtein, Tanja
Pasurek and Michael R\"{o}ckner are greatly appreciated. Part of this
research was carried out during my visits to the Department of Mathematics
of Bielefeld University. Financial support of these visits by the DFG
through SFB 701 \textquotedblleft Spektrale Strukturen und Topologische
Methoden in der Mathematik\textquotedblright\ and Alexander von Humboldt
Stiftung is gratefully acknowledged.

\end{document}